# ON THE POWERS OF SOME NEW CHI-SQUARE TYPE STATISTICS

**CLEMENT AMPADU**

Department of Mathematics
Central Michigan University
Mount Pleasant, MI 48859
U. S. A.
e-mail: ampad1cb@cmich.edu

**Abstract**

In this paper, four new Chi-Square type statistics are presented for testing the hypothesis of a uniform null versus specified trend alternatives. The powers of these test statistics are compared with the powers of the statistics considered by Steele and Chaseling [8]. The four test statistics are shown to have superior or equivalent powers to the powers of the test statistics considered by the authors for certain trend alternatives and for certain conditions placed on the cell count.

## 1. Introduction

Let $P_1 \ldots, P_k$ denote $k$ multinomial distribution cell probabilities. Let $P_i$ denote the fraction of discrete data falling into cell $i$. Let $O_1, \ldots, O_k$ denote the respective observed frequencies of the $k$ cells in a random sample of size $N$, then $\sum_{i=1}^{k} O_i = N.$ We are interested in testing the null hypothesis $H_0 : P_1 = \cdots = P_k$ versus several kinds of alternatives, in the case where a natural ordering of the cells exist. Let $E_i$ denote the expected cell count (expected frequency) of cell $i$ under $H_0$ for $i = 1, \ldots, k$.

Pearson [6], Neyman [5], Wilks [9] and Kullback [4] proposed the following Chi-Square type statistics for goodness-of-fit testing, respectively:

- $T_{X^2} = \sum_{i=1}^{k} \frac{(O_i - E_i)^2}{E_i}$ (Pearson),

- $T_{X^2(M)} = \sum_{i=1}^{k} \frac{(O_i - E_i)^2}{O_i}$ (Neyman),

- $T_{G^2} = 2\sum_{i=1}^{k} O_i \ln\left(\frac{O_i}{E_i}\right)$ (Wilks),

- $T_{G^2(M)} = 2\sum_{i=1}^{k} E_i \ln\left(\frac{E_i}{O_i}\right)$ (Kullback).

Let $\overline{E} = \frac{1}{k}\sum_{i=1}^{k} E_i$ and $\overline{O} = \frac{1}{k}\sum_{i=1}^{k} O_i$. We consider the following versions of the statistics above in our power study:

- $A_{X^2} = \sum_{i=1}^{k} \frac{(O_i - E_i)^2}{\overline{E}}$ (New Pearson),

- $A_{X^2(M)} = \sum_{i=1}^{k} \frac{(O_i - E_i)^2}{\overline{O}}$ (New Neyman),

- $A_{G^2} = 2\sum_{i=1}^{k} O_i \ln\left(\frac{O_i}{\overline{E}}\right)$ (New Wilks),

- $A_{G^2(M)} = 2\sum_{i=1}^{k} \overline{E} \ln\left(\frac{\overline{E}}{O_i}\right)$ (New Kullback).

Now let $Z_i = \sum_{j=1}^{i}(O_j - E_j), H_i = \sum_{j=1}^{i} E_j$ and $\overline{Z} = \sum_{j=1}^{k} Z_j p_j$. In conjunction with the original Pearson's Chi-Square statistic above, Steele and Chaseling [8] also consider in their power study the following discrete versions of goodness-of-fit statistics for continuous distributions developed by Kolmogorov [3], Smirnov [7], Cramer [2] and Anderson and Darling [1]:

- $S = \max_{1 \le i \le k} | Z_i |$ (Discrete Kolmogorov-Smirnov),
- $W^2 = N^{-1} \sum_{i=1}^{k} Z_i^2 p_i$ (Discrete Cramer-von Mises),
- $U^2 = N^{-1} \sum_{i=1}^{k} (Z_i - \overline{Z})^2 p_i$ (Discrete Watson),
- $A^2 = N^{-1} \sum_{i=1}^{k} \frac{Z_i^2 p_i}{H_i(1 - H_i)}$ (Discrete Anderson-Darling),
- $NS = \frac{1}{2} \sum_{i=1}^{k} | O_i - E_i |$ (Nominal Kolmogorov-Smirnov).

## 2. The Power Study

To make the comparison of the new statistics with those used by Steele and Chaseling [8] uniform, the powers of each of the ten test statistics are also approximated for a uniform null across 10 cells against the selected alternatives used by them. We also use sample sizes of 10, 20, 30, 50, 100 and 200, which are generated by using frequencies under the uniform null distribution of 1, 2, 3, 5, 10, and 20 per cell, respectively.

As done by Steele and Chaseling [8], the power of each test statistic is estimated from 10,000 simulated random samples. As observed by the authors, the discrete nature of the null distribution of each test statistic means that a critical value and corresponding power at a significance level of exactly 0.05 may not be possible. To overcome this, as done by the authors, the powers are obtained for critical values on both sides of the 0.05 level, and linear interpolation between them is used to derive an approximate power for the 0.05 level.

The table below summarizes the distributions as used by Steele and Chaseling [8] for the power study.

| Description | 1 | 2 | 3 | 4 | 5 | 6 | 7 | 8 | 9 | 10 |
|---|---|---|---|---|---|---|---|---|---|---|
| Uniform | 0.10 | 0.10 | 0.10 | 0.10 | 0.10 | 0.10 | 0.10 | 0.10 | 0.10 | 0.10 |
| Decreasing | 0.32 | 0.13 | 0.10 | 0.08 | 0.07 | 0.07 | 0.06 | 0.06 | 0.05 | 0.05 |
| Step | 0.05 | 0.05 | 0.05 | 0.05 | 0.05 | 0.15 | 0.15 | 0.15 | 0.15 | 0.15 |
| Tringular ($v$) | 0.17 | 0.13 | 0.10 | 0.07 | 0.03 | 0.03 | 0.07 | 0.10 | 0.13 | 0.17 |
| Platykurtic | 0.04 | 0.11 | 0.11 | 0.12 | 0.12 | 0.12 | 0.12 | 0.11 | 0.11 | 0.04 |
| Leptokurtic | 0.05 | 0.05 | 0.05 | 0.05 | 0.30 | 0.05 | 0.05 | 0.05 | 0.05 | 0.05 |
| Bimodal | 0.05 | 0.11 | 0.17 | 0.11 | 0.06 | 0.06 | 0.11 | 0.17 | 0.11 | 0.05 |

For the remainder of this paper we shall refer to the Discrete Kolmogorov-Smirnov statistic as Discrete KS, the Discrete Cramer-von Mises statistic as Discrete CVM, the Discrete Anderson-Darling statistic as Discrete AD, the Nominal Kolmogorov-Smirnov statistic as Nominal KS and the original Pearson's Chi-Square statistic as Chi-Square.

## 3. Powers of the Test Statistic for Each Alternative

### 3.1. Decreasing trend alternative distribution

With the decreasing trend alternative, Figure 1 indicates that the powers of Discrete AD rank the highest. The powers of Discrete CVM and Discrete KS rank second and third, respectively. The powers of New Pearson and the Chi-Square statistics are similar; we give them the same rank, four. The powers of Discrete Watson, New Wilks, Nominal KS and New Neyman, rank fifth, sixth, seventh and eighth, respectively.

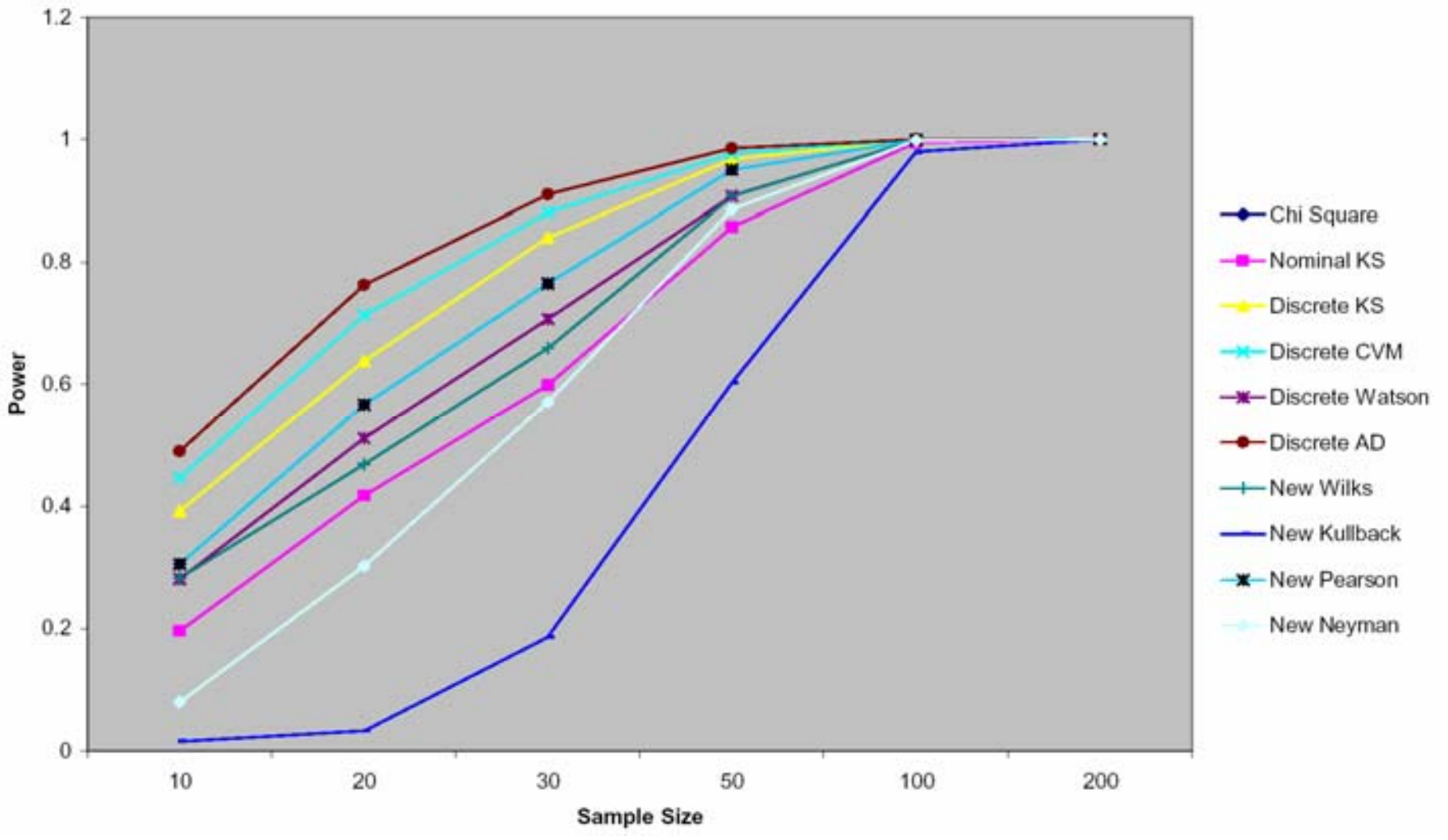

**Figure 1.** Power of the test statistics for a uniform null and decreasing alternative.

Comparing the new statistics with those used by Steele and Chaseling [8], we find the following: (1) The New Pearson statistic performs better than the Discrete Watson and the Nominal KS statistics. (2) The New Wilks statistic performs better than the Nominal KS statistics. (3) The New Neyman and the New Kullback statistics do not perform better than any of the statistics used by the authors. (4) In general the powers of all the new statistics are very high for sample sizes of at least five observations per cell $(N \geq 50)$ under the uniform null distribution, and they are identical to the powers of the test statistics used by the authors. (5) For sample sizes of at most three observations per cell $(N \leq 30)$ the New Pearson statistic gives the best power and has comparable power among the highest powers given by the test statistics used by the authors.

### 3.2. Step type alternative distribution

With the step-type trend alternative, Figure 2 indicates that the powers of Discrete CVM and Discrete KS are the highest and are similar. The powers of Discrete AD are the second highest. For sample sizes of approximately less than 3 per cell) 30 $(N < 30)$, the powers of New

Neyman are the third highest, and for sample sizes of at least approximately 3 per cell $(N \geq 30)$, the powers of Discrete Watson are the third highest. The powers of New Pearson, New Wilks, Nominal KS, and the Chi-Square statistics are the fourth highest, and are all similar. The powers of the New Kullback statistic lag behind the rest.

Comparing the new statistics with those used by Steele and Chaseling [8] we find the following: (1) If the sample size is less than 3 observations per cell $(N < 30)$, New Neyman performs better than the Chi-Square and the Nominal KS statistics used by the authors. (2) New Pearson and New Wilks appear to have the same statistical power compared to that of the Nominal KS and Chi-Square statistics used by the authors. (3) New Kullback does not perform better than any of the statistics used by the authors. (4) Among all the new statistics, the New Neyman statistic is more powerful in detecting the step type trend (irrespective of the sample size, it has the highest power among all the new statistics). (5) For samples sizes of at least 10 per cell $(N \geq 100)$ the powers of all the new statistics are very high under the uniform null, and are identical to the statistics used by the authors.

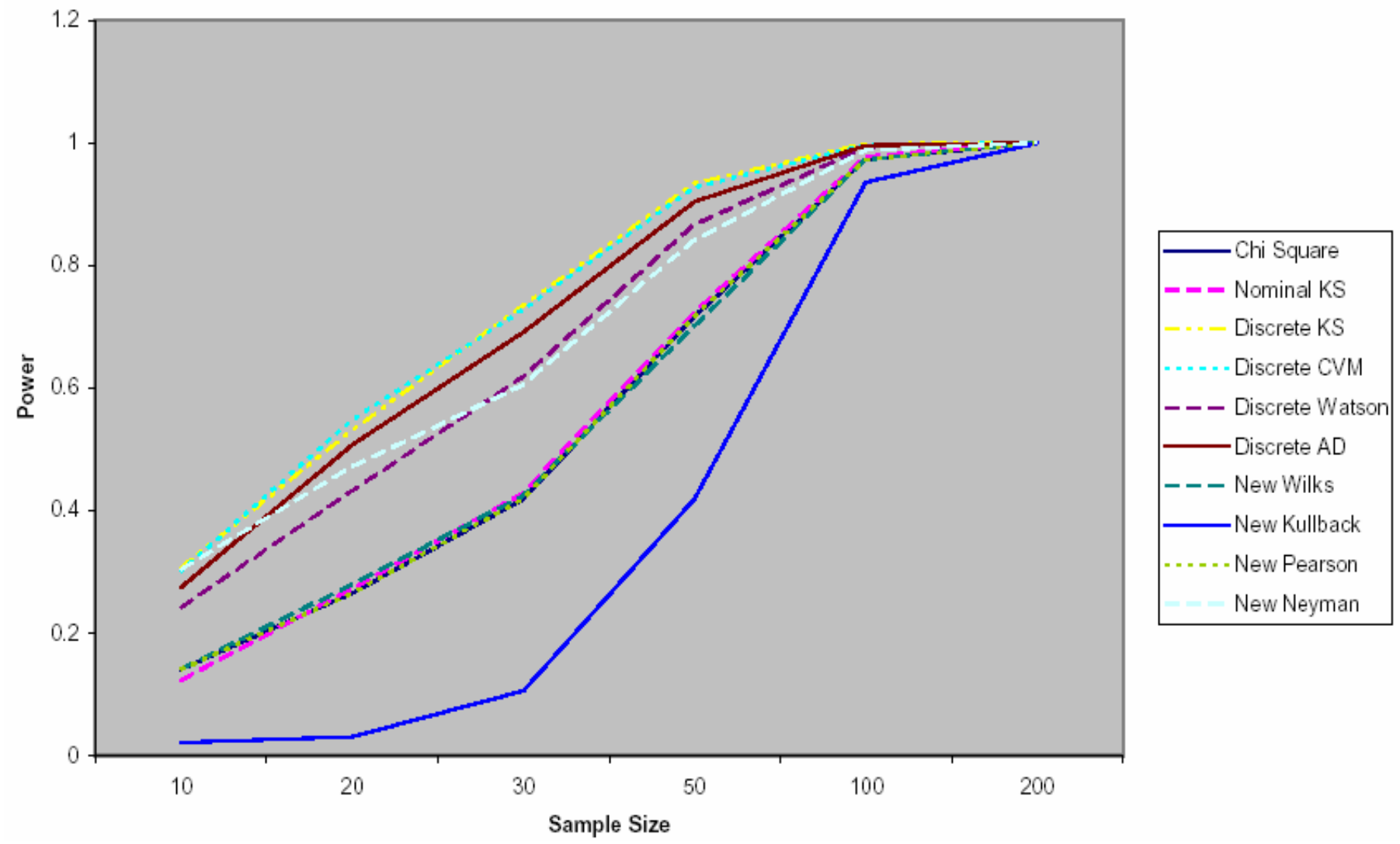

**Figure 2.** Power of the test statistics for a uniform null and step-type alternative.

### 3.3. Triangular alternative distribution

With the triangular ($v$) trend alternative, Figure 3 shows that the powers of New Pearson and Chi-Square are the highest, and are similar. The New Wilks statistic has the second highest power, New Neyman has the third highest power, Nominal KS has the fourth highest power, Discrete Watson has the fifth highest power, Discrete KS has the sixth highest power, Discrete CVM has the seventh highest power, Discrete AD has the eighth highest power, and New Kullback lags behind the rest.

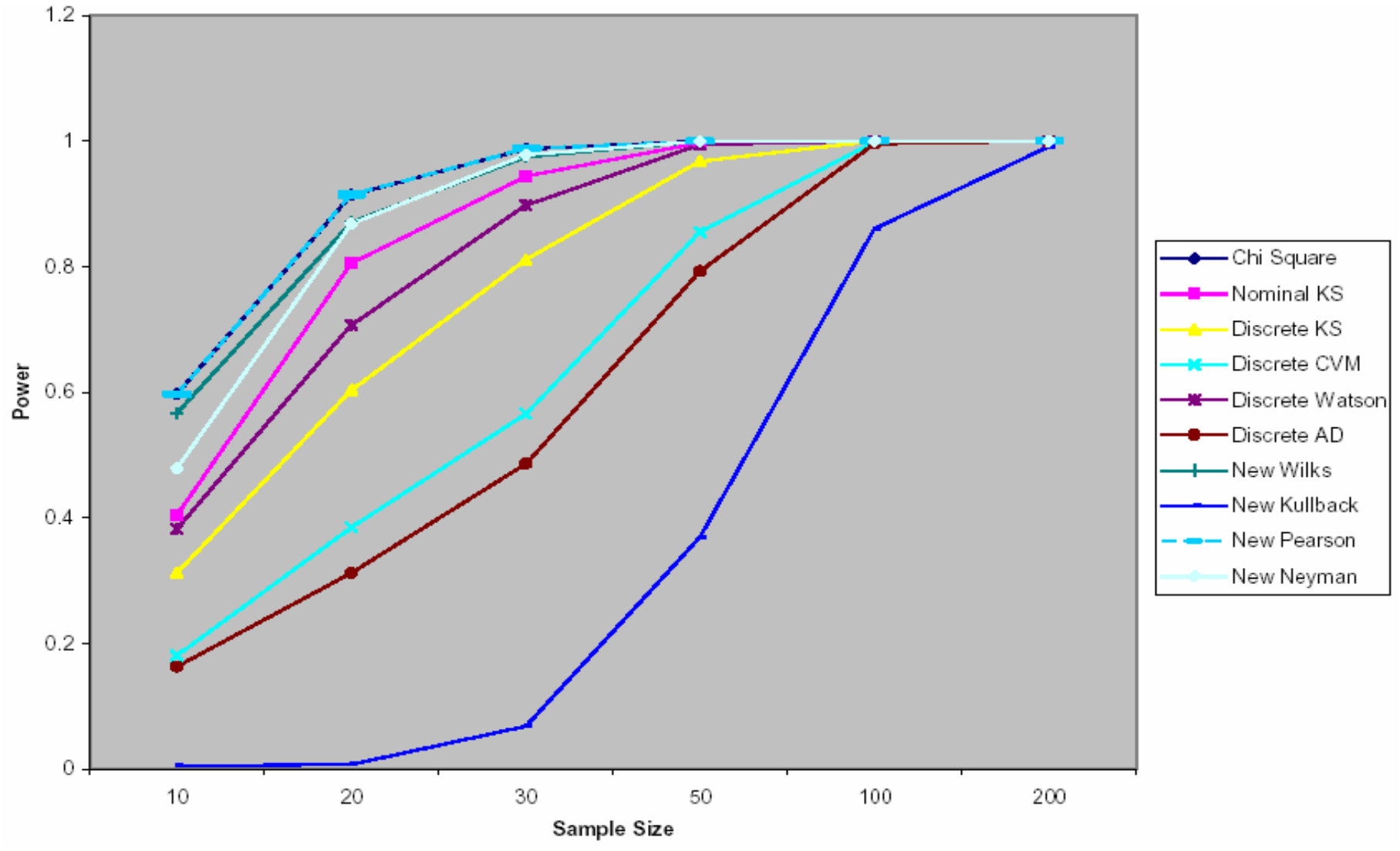

**Figure 3.** Power of the test statistics for a uniform null and triangular ($v$) alternative.

Comparing the new statistics with those used by Steele and Chaseling [8], we find the following: (1) The New Pearson, New Neyman, and New Wilks statistics have superior power than those considered by the authors, with the exception of Chi-Square whose power is comparable to New Pearson. (2) For samples sizes of at least 10 per cell ($N \geq 100$) the powers of all the new statistics are very high under the uniform null, and with the exception of New Kullback are identical to all the statistics used by the authors. (4) Among all the new statistics, New Pearson detects the triangular trend best, it has the highest power.

### 3.4. Platykurtic alternative distribution

With the platykurtic trend alternative, Figure 4 indicates the following: (1) If the sample size is less than 3 per cell, the powers of Discrete Watson rank the highest. The powers of New Wilks and Nominal KS rank second and third, respectively. The powers of New Pearson and Chi Square are similar, and are the fourth highest. The powers of New Neyman rank the fifth highest. The powers of Discrete AD, Discrete CVM, Discrete KS and New Kullback are similar, and rank sixth. (2) If the sample size is at least 3 per cell but less than 5 per cell, the powers of Discrete Watson rank the highest. The powers of New Wilks are the second highest. The powers of Nominal KS, New Pearson and Chi-Square are similar, and are the third highest. The powers of New Neyman are the fourth highest, whilst those of New Kullback are the fifth highest. The powers of Discrete AD, Discrete KS and Discrete CVM, rank sixth, seventh and eight, respectively. (3) If the sample size is at least 5 per cell and less than 10 per cell, New Wilks and Discrete Watson have the highest power and the powers of New Kullback rank second. The powers of New Pearson and Chi-Square are similar; we give them the same rank, three. The powers of New Neyman, Nominal KS, Discrete AD, Discrete CVM and Discrete KS have rank four, five, six, seven and eight, respectively. (4) If the sample size is at least 10 per cell, the powers of New Kullback, New Wilks and Discrete Watson are the highest. The powers of New Neyman, New Pearson and Chi-Square are similar and rank the second highest. The powers of Nominal KS, Discrete AD, Discrete CVM and Discrete KS have rank three, four, five and six, respectively.

Comparing the new statistics with those used by Steele and Chaseling [8] we find the following: (1) If the sample size is less than 3 per cell, New Wilks, New Pearson and New Neyman perform better than the Discrete AD, Discrete CVM and Discrete KS statistics used by the authors. New Wilks performs better than the Nominal KS and Chi-Square statistics used by the author. Among all the new statistics, New Wilks has the highest power and appears to have a better chance of detecting the platykurtic trend. (2) If the sample size is at least 3 per cell but less than 5 per cell, all the new statistics perform better than the Discrete AD, Discrete KS and Discrete CVM statistics used by the authors. New Wilks performs better than the Nominal KS and Chi-Square statistics used by the authors. Among, all the new statistics, New Wilks has the highest power and appears to have a better chance of detecting the platykurtic trend. (3) If the sample size is at least 5 per cell but less than 10 per cell, all the new statistics perform better than the Nominal KS, Discrete AD, Discrete CVM and Discrete KS statistics used by the authors. New Wilks and New Kullback perform better than Chi-Square. Among all the statistics (not just the new ones), New Wilks has one of the best chances of detecting the platykurtic trend, it has rank one in conjunction with the Discrete Watson statistic used by the authors. (4) If the sample size is at least 10 per cell, all the new statistics perform better than the Nominal KS, Discrete AD, Discrete CVM and Discrete KS statistics considered by the authors. New Kullback and New Wilks perform better than the Chi-Square statistics used by the authors. Among all the statistics (not just the new ones), New Wilks and New Kullback have one of the best chances of detecting the platykurtic trend, these statistics have rank one in conjunction with the Discrete Watson statistics used by the authors.

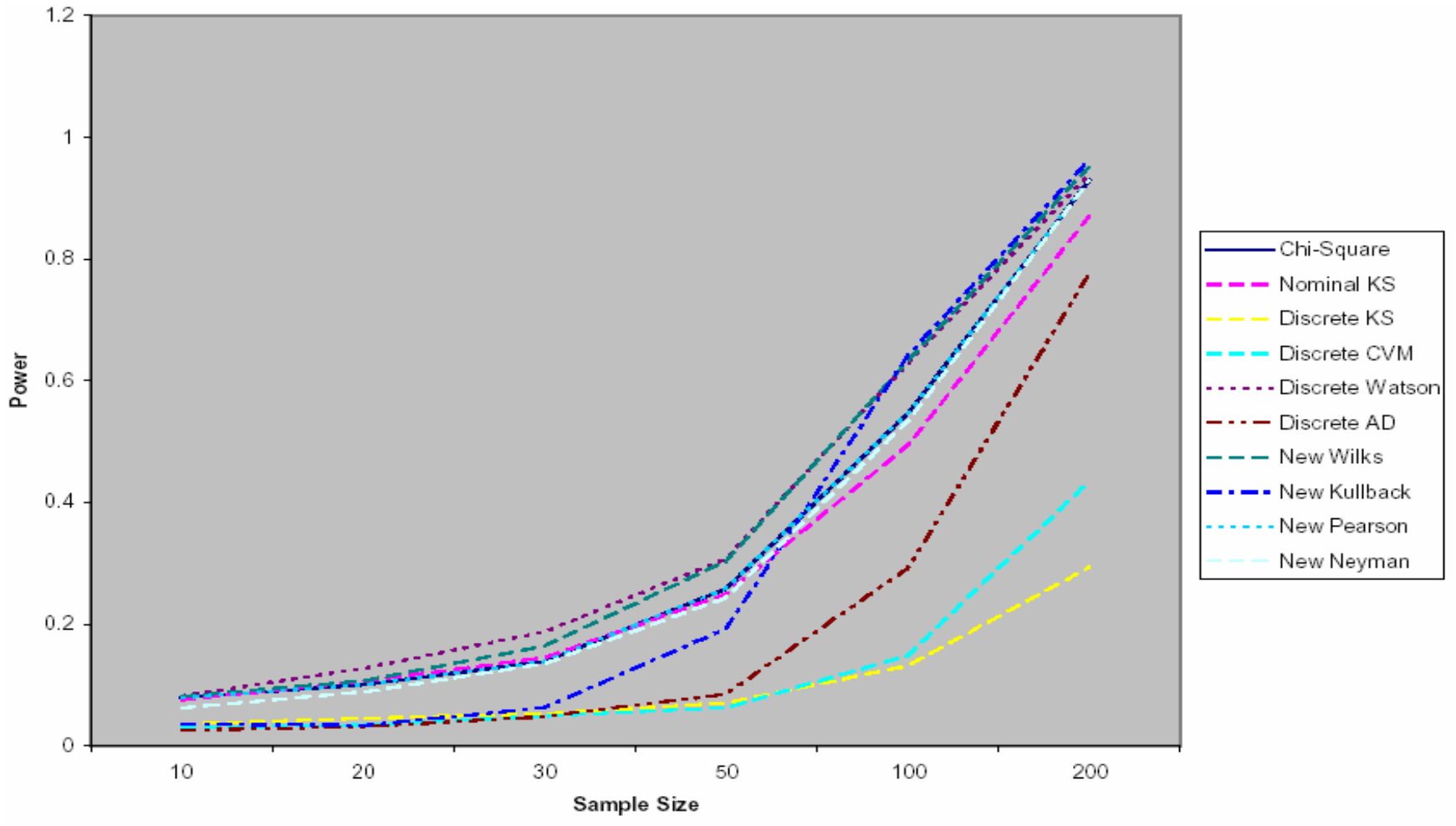

**Figure 4.** Powers of the test statistics for a uniform null and platykurtic alternative.

**3.5. Leptokurtic alternative distribution**

With the leptokurtic trend alternative, Figure 5 shows that the powers of Discrete Watson rank the highest. The powers of New Pearson and Chi-Square are the second highest, and are similar. For sample sizes of approximately less than 2 per cell $(N < 20)$, the powers of New Wilks rank the third highest, and for sample sizes of approximately at least 2 per cell $(N \geq 20)$, the powers of New Neyman rank the third highest. For sample sizes of approximately less than 2 per cell $(N < 20)$, the powers of New Neyman rank the fourth highest. For sample sizes of approximately less than 3 per cell $(N < 30)$, the powers of Nominal KS rank the fifth highest, and for sample sizes of at least 3 per cell $(N \geq 30)$, the powers of New Wilks and Nominal KS rank the fifth highest, and are similar. The powers of Discrete KS, Discrete CVM and Discrete AD, rank sixth, seventh and eight, respectively. The powers of the New Kullback statistic lag behind the rest.

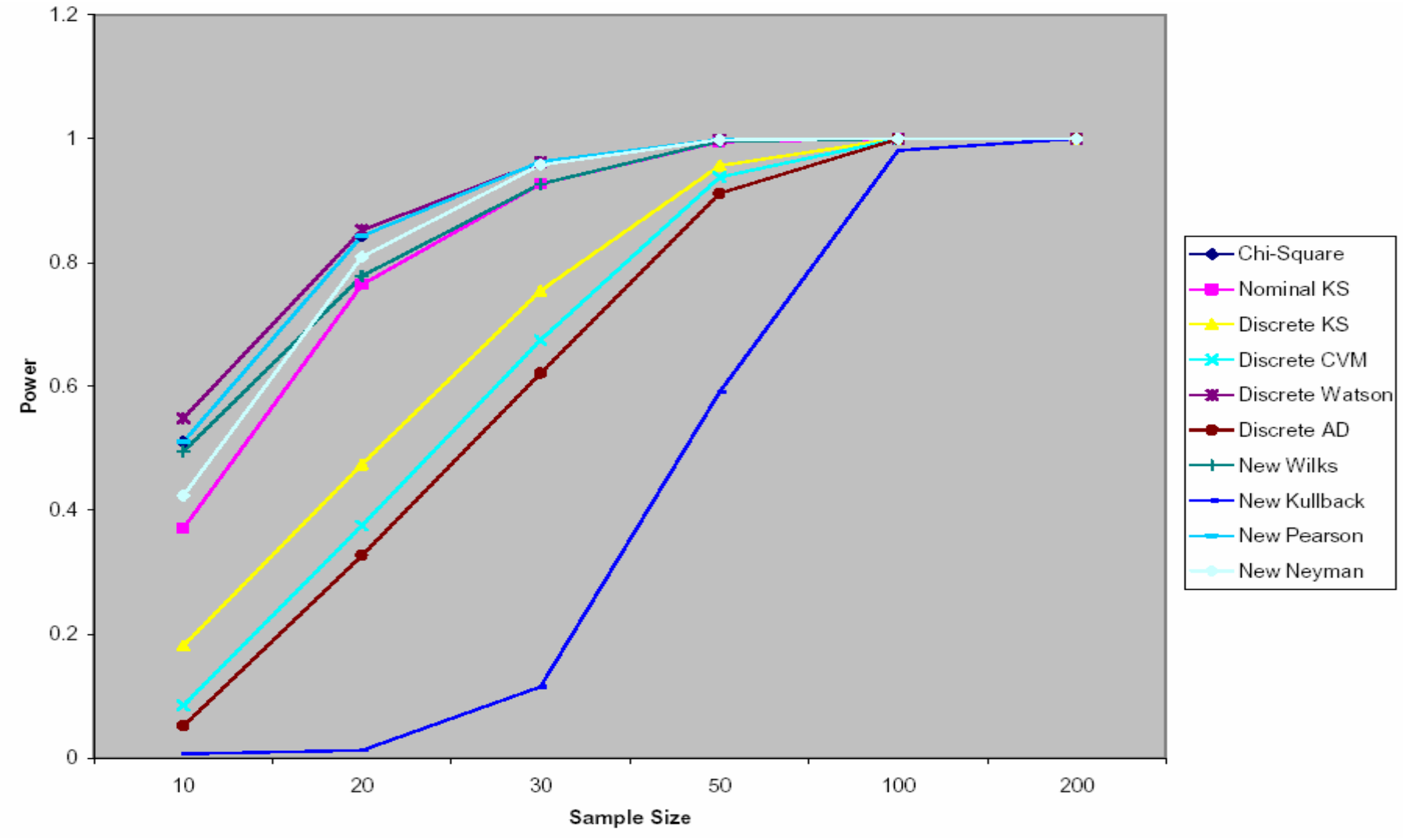

**Figure 5.** Powers of the test statistic for a uniform null and leptokurtic alternative.

Comparing the new statistics with those used by Steele and Chaseling [8], we find the following: (1) New Kullback does not perform better than any of the statistics used by the authors. (2) New Pearson,

New Wilks and New Neyman have superior power compared to the Discrete KS, Discrete CVM, and the Discrete AD statistics used by the authors; the only exception is with the Discrete Watson and Chi-Square statistics used by the authors. (3) For sample sizes of at least 10 per cell $(N \geq 100)$, the powers of all the new statistics are very high under the uniform null, and are identical to the statistics used by the authors. (4) Among all the new statistics, New Pearson detects the leptokurtic trend best, it has the highest power.

### 3.6. Bimodal alternative distribution

With the bimodal trend alternative, Figure 6 indicates the following: (1) For sample sizes less than 3 per cell $(N < 30)$, New Pearson, New Neyman, Nominal KS, New Wilks and Chi-Square have the highest powers, and the powers are similar, we give these statistics the highest ranking, one. Discrete Watson and New Kullback have the second and third highest powers, respectively. Discrete KS, Discrete CVM and Discrete AD have similar powers and we rank them fourth. (2) For sample sizes of at least 3 per cell and less than 5 per cell $(30 \leq N < 50)$, New Pearson, New Neyman, Nominal KS, New Wilks and Chi-Square have the highest powers, and the powers are similar, we give these statistics the highest ranking, one. New Kullback and Discrete Watson are ranked second and third, respectively. Discrete KS is ranked fourth, and Discrete CVM and Discrete AD have similar statistical power, and are ranked fifth. (3) For sample sizes of at least 5 per cell and less than 10 per cell $(50 \leq N < 100)$, New Pearson, New Neyman, Nominal KS, New Wilks and Chi-Square have the highest powers, and the powers are similar, we give these statistics the highest ranking, one. New Kullback and Discrete Watson are ranked second and third, respectively. Discrete KS is ranked fourth, Discrete AD is ranked fifth, and Discrete CVM has the last rank, sixth. (4) For sample sizes of at least 10 per cell $(N \geq 100)$, Chi-Square, New Pearson, New Neyman, Nominal KS and New Wilks have the highest powers which are similar, we give these statistics the highest ranking, one. New Kullback, Discrete Watson, Discrete AD, Discrete KS and Discrete CVM are ranked second, third, fourth, fifth and sixth, respectively.

Comparing the new statistics with those used by Steele and Chaseling [8], we find the following: (1) For sample sizes less than 3 per cell, all the new statistics perform better than the Discrete KS, Discrete CVM and Discrete AD statistics used by the authors. New Pearson, New Neyman and New Wilks perform better than the Discrete Watson statistic used by the authors. With the exception of New Kullback, all the other new statistics have rank one and thus have the best chance among them to detect the bimodal trend alternative. (2) For sample sizes of at least three per cell and less than 5 per cell, all the new statistics perform better than the Discrete Watson, Discrete KS, Discrete CVM and Discrete AD statistics used by the authors. With the exception of New Kullback, all the other statistics have rank one, and thus have the best chance of detecting the bimodal trend alternative. (3) For sample sizes of at least 5 per cell and less than 10 per cell, all the new statistics perform better than the Discrete Watson, Discrete KS, Discrete AD and Discrete CVM statistics used by the authors. With the exception of New Kullback all the new statistics have rank one, and thus have the best chance of detecting the bimodal trend alternative. (4) For sample sizes of at least 10 per cell, all the new statistics perform better than the Discrete Watson, Discrete AD, Discrete KS, and Discrete CVM statistics used by the authors. With the exception of New Kullback, all the new statistics have rank one, and thus have the best chance of detecting the bimodal trend alternative.

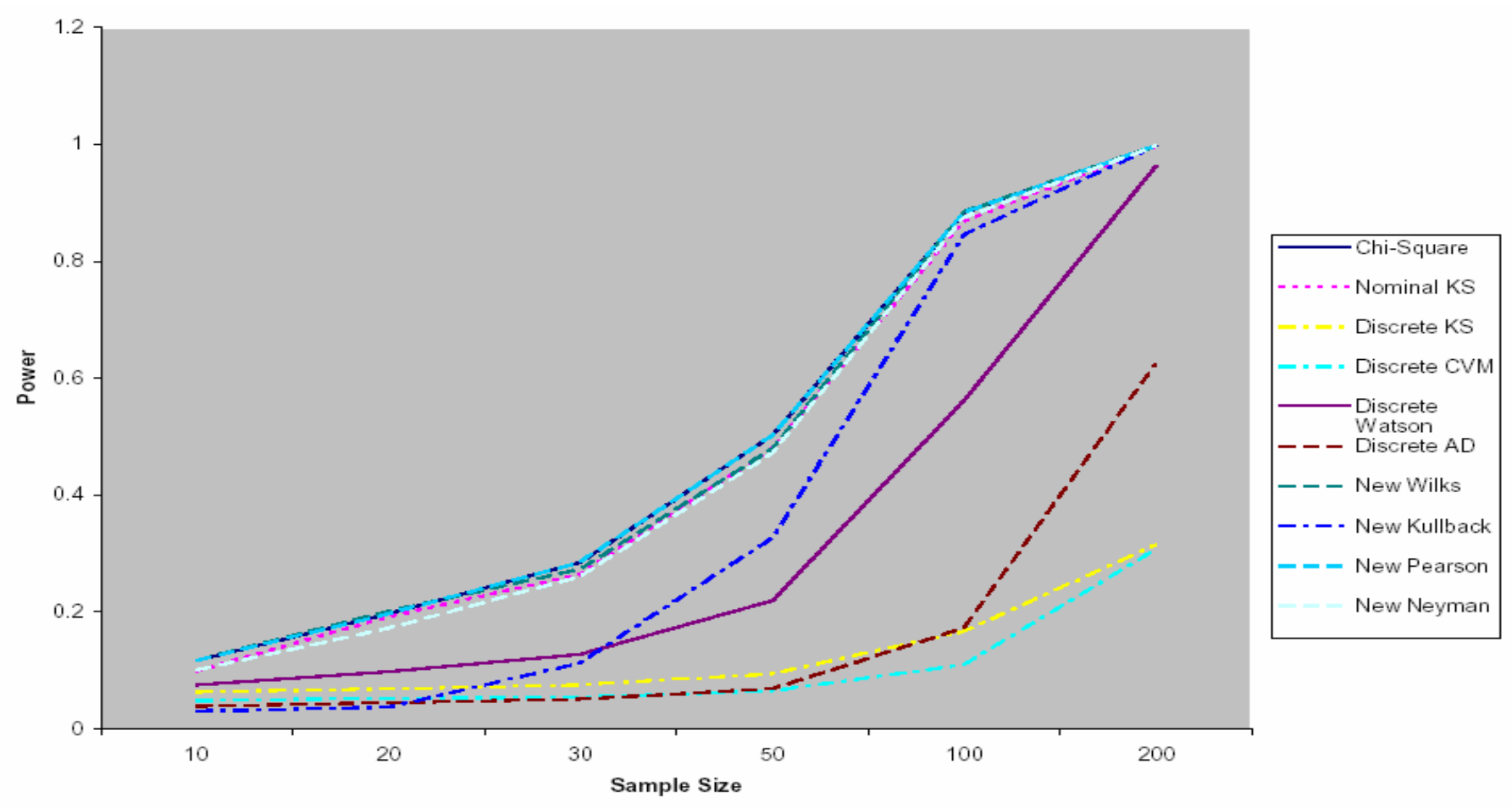

**Figure 6.** Powers of the test statistic for a uniform null and the bimodal alternative.

## 4. Conclusions

Under the same conditions as the power study by Steele and Chaseling [8], we find that for certain trend alternatives and also for certain conditions placed on the cell counts, that the new test statistics will optimize the power of the test of goodness of fit much better than the statistics used by the authors, or will provide the same optimum power as the statistics used by the authors.

From the power analysis in Section 3, we see that the situations under which these happen are as follows: (1) When the alternate trend is triangular, New Pearson detects the trend just as effectively as Chi-Square. (2) When the trend is platykurtic and the cell count is between 5 and 10 per cell, New Wilks detects the trend just as effectively as Discrete Watson. (3) When the trend is platykurtic and the cell count is at least 10 per cell, New Kullback and New Wilks detect the trend just as effectively as Discrete Watson. (4) When the trend is Bimodal, New Pearson, New Neyman and New Wilks detects the trend much better than any of the statistics considered by Steele and Chaseling [8]. However, since the bimodal trend deviates from the uniform null, it is unclear whether the superior performance of the new statistics over those considered by Steele and Chaseling [8] is attributable to the deviation or to the specific shape.

| | |
|---|---|
| Kindly return the proof after correction to:<br><br>*The Publication Manager*<br>*Far East J. Theo. Stat.*<br>*Pushpa Publishing House*<br>*Vijaya Niwas*<br>198, *Mumfordganj*<br>*Allahabad*-211002 (*India*)<br><br>along with the print charges* by the <u>fastest mail</u><br><br>***Invoice attached** | Proof read by: .......................................<br>Signature: ..........................................<br>Date: ................................................<br>Tel: ..................................................<br>Fax: .................................................<br>E-mail: .............................................<br>Number of additional reprints required ................................<br>Cost of a set of 25 copies of additional reprints @ EURO 12.00 per page.<br>(25 copies of reprints are provided to the corresponding author ex-gratis) |